\documentclass{conm-p-l}

\usepackage{amssymb}

\usepackage{}

\newtheorem{theorem}{Theorem}[section]
\newtheorem{lemma}[theorem]{Lemma}

\theoremstyle{definition}
\newtheorem{definition}[theorem]{Definition}
\newtheorem{example}[theorem]{Example}

\theoremstyle{remark}
\newtheorem{remark}[theorem]{Remark}

\numberwithin{equation}{section}

\begin{document}

\title{$H_\infty\neq E_\infty$}


\author{Justin Noel}
\address{Max Planck Institute for Mathematics, Vivatsgasse 7, Bonn, Germany}
\email{justin@mpim-bonn.mpg.de}
\thanks{The author was partially supported by the Max Planck Institute for Mathematics and by the Deutsche Forschungsgemeinschaft through Graduiertenkolleg 1150.}

\subjclass[2010]{55P43, 55P47, 55P48}

\date{\today}

\begin{abstract}
We provide an example of a spectrum over $S^0$ with an $H_\infty$
structure which does not rigidify to an $E_3$ structure.  It follows
that in the category of spectra over $S^0$ not every $H_\infty$ ring
spectrum comes from an underlying $E_\infty$ ring spectrum.  After
comparing definitions, we obtain this example by applying
$\Sigma^\infty_+$ to the counterexample to the transfer conjecture
constructed by Kraines and Lada.
\end{abstract}

\maketitle

\section{Introduction}

In recent years there has been a renewed interest in the study of
$E_\infty$ ring spectra and their strictly commutative analogues,
commutative $S$-algebras. Such spectra are equipped with a well-behaved
theory of power operations.  This structure provides formidable
computational tools which can be used to deduce a number of surprising
results (for some examples see \cite[Ch.~2]{BMMS86}).

Such operations determine and are determined by an $H_\infty$ ring 
structure, the analogue of an $E_\infty$ ring structure in the
stable \emph{homotopy} category.  The theory of power operations is
sufficiently rich that one might conjecture that  every $H_\infty$ ring
spectrum is obtained by taking an $E_\infty$ ring spectrum and then
passing to the homotopy category.  

This turns out to be a stable analogue of the transfer conjecture, a
conjectural equivalence between the homotopy category of infinite loop
spaces and a subcategory of the homotopy category of based spaces whose
objects admit certain transfer homomorphisms (see \cite{KrL79} for a
more complete description).   

Kraines and Lada demonstrate the falsehood of the transfer conjecture by
contructing an explicit counterexample. In their paper, Kraines and Lada
define the notion of an $L(n)$ space. When $n=2,$ this is a space
equipped with transfer homomorphisms. They make use of the
following implications 
\begin{align*}
   X\text{ is an infinite loop space} &\implies X\text{ is an }E_\infty\text{ space}\\
  &\iff X\text{ is an } L(\infty)\text{ space}\\
  &\implies X\text{ is an }L(n)\text{ space}\\ 
  &\implies X\text{ is an }L(n-1)\text{ space}\ldots\\
\end{align*}
We can view lifting an $L(n)$ structure to an $L(n+1)$ structure and so on as constructing an action of the $E_\infty$ operad up to increasingly coherent homotopy.

\begin{theorem}[{\cite{KrL79}}]\label{thm:KL} Let $s$ be a generator of
  $\mathrm{Prim}H^{30}(BU;\mathbb{Z}_{(2)}).$  Define $KL$ by the
  following fibration sequence:  
  \[KL\xrightarrow{i} BU_{(2)} \xrightarrow{4s} K(\mathbb{Z}_{(2)},30). \] 
   Then $i$ is a map of $L(2)$ spaces, but the $L(2)$ structure on $KL$
   does not lift to an $E_3$ structure. In particular, $KL$ does not
   admit an $E_\infty$ structure compatible with this $L(2)$ structure.
\end{theorem}

After some translation we will prove the following theorem, which
provides an example in the category of $H_\infty$ ring spectra augmented
over $S^0$ whose $H_\infty$ structure does not arise by
forgetting an $E_\infty$ structure. 

\pagebreak[2]
\begin{theorem}\label{thm:main}
   The map \[\Sigma^\infty_+KL \xrightarrow{\Sigma^\infty_+i}
   \Sigma^\infty_+ BU_{(2)}\]
   is a map of $H_\infty$ ring spectra augmented over $S^0$, but the
   $H_\infty$ ring structure on $\Sigma^\infty_+ KL$ does not lift to an
   $E_3$ structure. In particular, $\Sigma^\infty_+ KL$ does not admit a
   compatible $E_\infty$ ring structure.
\end{theorem}

To prove this we will show that $\Sigma^\infty_+$ takes $L(2)$ spaces to
$H_\infty$ ring spectra under $S^0$ and takes $E_\infty$ spaces
spaces to $E_\infty$ ring spectra under $S^0$.  This comparison is
deduced immediately from some of the results in \cite{May09}.  

The author would also like to thank Peter May and the anonymous referee for
their helpful comments and suggestions concerning this paper.

\section{$L(n)$ spaces and spectra}
Let $\mathcal{L}$ be the linear isometries operad. We will abuse
notation and let $L$ denote the associated reduced monad on pointed
spaces with Cartesian products, spaces under $S^0$ with smash products,
and spectra under $S^0$ with smash products.

In particular:
\begin{itemize}
\item $L$ is an endofunctor on pointed spaces satisfying
  \[LY=\coprod_{n\geq0}\mathcal{L}(n)\times_{\Sigma_n}Y^n/(\sim),\] where
  $\sim$ represents the obvious base point identifications.
\item $L$ is an endofunctor on spaces under $S^0$ satisfying
  \[LY=\coprod_{n\geq0}\mathcal{L}(n)\times_{\Sigma_n}Y^n/(\sim),\] where
  $\sim$ represents the obvious unit map identifications.
\item $L$ is an endofunctor on the Lewis-May-Steinberger category of
  spectra (see \cite{LMS86}) under $S^0$ satisfying
  \[LE=\bigvee_{n\geq0}\mathcal{L}(n)\ltimes_{\Sigma_n}E^{\wedge
  n}/(\sim),\] where $\sim$ represents the obvious unit map
  identifications (see \cite[4.9,6.1]{EKMM97}).
\end{itemize}

We justify this abuse of notation with the following lemma:
\begin{lemma}[{\cite[4.8, p.~1027]{May09}}]\label{lem:QD}
  We have the following chain of isomorphisms natural in based
  spaces\footnote{It is helpful to think of this basepoint as the  
  multiplicative unit.}
  $X$ \begin{align*}
   \Sigma^\infty_+LX &\equiv \Sigma^\infty(LX)_+\\
   	&\cong \Sigma^\infty L(X_+)\\
   	&\cong L\Sigma^\infty X_+\\
   	&\equiv L\Sigma^\infty_+ X.\\
 \end{align*}
\end{lemma}

For simplicity, for the remainder of this paper we will assume all
spaces are non-degenerately based and let $e\colon \textrm{Id}\rightarrow L$ and $\mu\colon L^2=LL\rightarrow L$ denote the structure maps of $L$.

Recall that the category of $L$-algebras in group-like pointed spaces is
equivalent to the category of infinite loop spaces. The following
definition provides a categorical filtration between spaces and
homotopy coherent $L$-algebras (which are weakly equivalent to
$L$-algebras).

\begin{definition}
  A based space $X$ is $L(n)$ if one can construct maps 
  \[ f_k\colon I^k\times L^{k+1}X\rightarrow X \textrm{ for } k<n\] such that
  \begin{enumerate}
  \item the composite $X\xrightarrow{e} LX\xrightarrow{f_0} X$ is the identity,
  \item if $t_j=0$, $f_k(t_1,\ldots, t_k, z)=f_{k-1}\circ(\mathrm{Id}_{I^{k-1}}\times L^{j-1}\mu L^{k-j}) (t_1,\ldots,\widehat{t_j},\ldots, t_k,z)$,
  \item and if $t_j=1$, $f_k(t_1,\ldots, t_k, z)=f_{j-1}\circ(\mathrm{Id}_{I^{j-1}}\times L^{j}f_{k-j}) (t_1,\ldots,\widehat{t_j},\ldots, t_k,z)$.
  \end{enumerate}
\end{definition}

\begin{remark}
  Despite the similarity in notation, we remind the reader that the
  \emph{property} of being $L(n)$ has nothing to do with the
  \emph{space} $\mathcal{L}(n)$. We also note that being $E_n$ does not imply the space is $L(n)$.
\end{remark}

\begin{remark}
  Note that our definition of a $L(n)$ space is different from that
  of a $Q_n$ space used in \cite{KrL79}. Kraines and Lada restrict to
  the case when $X$ is connected, in which case $L$ could be replaced
  with $Q=\Omega^\infty\Sigma^\infty.$  In this respect, our definition
  is more general.
\end{remark}

We illustrate our definition with a sequence of examples (for more
detailed exposition and proofs see \cite{KrL79} or \cite[V]{CLM76}).
\begin{example}\ \nopagebreak
  \begin{enumerate}
	\item By definition, every based space is a $L(0)$ space.
	\item A based space $X$ is $L(1),$ if the canonical map
	  $X\rightarrow LX$ admits a retraction $\mu_X$, which we can regard
	  as the multiplication on $X$.
	\item \label{enum:h_2} A space $X$ is $L(2),$ if it is $L(1)$ and we
	  have a specified homotopy $I\times L^2X\rightarrow X$ between
	  $\mu_X\mu$ and $\mu(\mu).$ In other words, $X$ is an
	  $L$-algebra in the homotopy category of pointed spaces.
	\item \label{enum:h_3} In the case $L$ is the monad associated to an
	operad, by a result of Lada \cite{CLM76}, $X$ is $L(\infty)$ if and
	only if it admits an $L$-algebra structure. If the components of $X$
	form a group under the induced multiplication, then $X$ is $L(\infty)$
	if and only if it has the homotopy type of an infinite loop space.
  \end{enumerate}
\end{example}

There is an obvious analogue of the above definition with based spaces
replaced by spectra under $S^0$, where $L$ is the monad
whose algebras are $E_\infty$ ring spectra in this category
\cite[6.2]{May09}. So we obtain an analogous categorical filtration between
spectra under $S^0$ and $E_\infty$ ring spectra. 

Applying this equivalence to the definition of $L(n)$ spectra, we see
that the definition of an $L(2)$ spectrum is precisely the definition of
a $H_\infty$ ring spectrum under $S^0$ \cite{BMMS86}. By
Lemma \ref{lem:QD} we see that applying $\Sigma^{\infty}_+$ to the map
\[ KL\rightarrow BU \]
of $L(2)$ spaces constructed by Kraines and Lada we obtain a map of
$H_\infty$ ring spectra augmented over $S^0$. 

To see that $\Sigma^{\infty}_+KL$ is not an $E_{\infty}$ ring spectrum
we apply the argument of \cite[\S 8]{KrL79}. There they show that 
the Postnikov system for $KL$ gives rise to a fibration sequence:
\[  KL_{\leq 29}\rightarrow KL_{\leq 28}\simeq BU_{\leq
28}\xrightarrow{\tau} K(\mathbb{Z}/(4\cdot 15!), 30). \] If $KL$ were an
infinite loop space, $KL_{\leq 29}$ would be as well and the
$k$-invariant would be an infinite loop map.  However they demonstrate
that $\tau$ can not be delooped twice to a multiplicative map and so the
above Postnikov fibration can not be delooped twice to a Postnikov
system of $A_\infty$ spaces. As a consequence $KL_{\leq 29}$ and $KL$
fail to be $E_3$ spaces and the corresponding 
suspension spectra fail to have induced $E_3$ structures.

\bibliographystyle{amsplain}

\bibliography{/home/justin/work/notes/biblio/biblio}
\end{document}